# Reduced Precision Checking to Detect Errors in Floating Point Arithmetic

Yaqi Zhang, Ralph Nathan, and Daniel J. Sorin
Duke University

**Abstract**

*In this paper, we use reduced precision checking (RPC) to detect errors in floating point arithmetic. Prior work explored RPC for addition and multiplication. In this work, we extend RPC to a complete floating point unit (FPU), including division and square root, and we present precise analyses of the errors undetectable with RPC that show bounds that are smaller than prior work. We implement RPC for a complete FPU in RTL and experimentally evaluate its error coverage and cost.*

## 1 Introduction

In this paper, we focus on detecting errors in floating point arithmetic performed by the floating point units (FPUs) in processor cores. We focus on all of the common floating point operations: addition/subtraction, multiplication, division, and square root.

We base our error detection on the previously introduced idea of reduced precision checking (RPC) [1][9]. The key idea behind RPC is to check a full-precision (32-bit, in this paper) FPU with a reduced-precision checker FPU. The full-precision FPU operates on IEEE standard operands that have a sign bit, 8 exponent bits, and 23 fraction bits. The checker FPU operates on operands that also have a sign bit and 8 exponent bits, but the fractions have *k* bits, where *k*≤23 is a parameter that is chosen by the designer. In the absence of errors, the full-precision result of the FPU will match the reduced-precision result of the checker; only errors in the least significant bits can evade detection. RPC has an inherent trade-off between the cost of the checker and the likelihood of it not detecting an error.

Prior work has developed RPC for addition/subtraction [1] and multiplication [9]. The prior work introduced the idea of RPC and built and evaluated hardware implementations for addition/subtraction and multiplication. Prior work is limited, however, in several important ways. First, it does not include floating point division or square root, both of which are common hardware units. Second, we show that the size of the undetected errors for addition/subtraction and multiplication can be greatly reduced; prior work made overly conservative assumptions that led to unnecessarily large errors evading detection. Third, for some floating point input operations on some standard operands[1] (e.g., subtraction of two similarly sized numbers), the prior RPC schemes were unable to detect any errors at all and had to be disabled for these operation/operand combinations.

In this paper, we make the following contributions:

- We develop and evaluate RPC for a *complete FPU* in RTL—with addition/subtraction, multiplication, division, and square root. The results include error detection coverage, chip area, and energy.
- We mathematically compute the maximum undetectable error for every floating point operation, and we show that we can obtain tighter bounds for addition/subtraction and multiplication than prior work.
- Using a new RPC technique, that we call "reverse checking," we eliminate the unchecked operation/operand pairs of prior work. That is, our RPC unit checks every operation on all standard operands.

## 2 Background and Notation

In this section, we introduce background and notation that we use throughout the rest of this paper. We assume that all full-precision floating point numbers adhere to the IEEE-754 standard. Furthermore, we consider only 32-bit floating point numbers in this paper, but our approach applies in an identical way to 64-bit floating point numbers.

### 2.1 Floating Point Format

As illustrated in Figure 1, a 32-bit floating point number *X* has 1 sign bit (X[31]), 8 exponent bits (X[30:23]), and 23 fraction bits (X[22:0]). For notational clarity, we use the symbol ∃ to denote the **exponent bits** (i.e. X[30:23]). The sign bit with value of 0 or 1 indicates a positive or negative **sign ($S_X$)**, respectively, i.e., $S_X = (-1)^{X[31]}$. The base-2 exponent is represented in unsigned, biased notation, such that the **exponent of $X$ ($E_X$)** is the unsigned value of the exponent bits after 127 is subtracted (i.e., $E_X = \exists - 127$). The mantissa of *X* is an implicit 1 followed by the 23 fraction bits.[2] We define **fraction ($f_X$)** as the *integer value* of the lower 23 bits of *X* (X[22:0]), and **mantissa ($M_X$)** as the *floating point value* of the mantissa, including the implicit one, i.e. $1.\{f_X\}$. Curly braces denote concatenation. Thus the floating point value of *X* is:

$$(-1)^{X[31]} \times 2^{\exists-127} \times 1.\{X[22:0]\}$$
$$= S_X \times 2^{E_X} \times 1.\{f_X\}$$
$$= S_X \times 2^{E_X} \times M_X$$

[1] A non-standard operand is a denorm, infinity, or NaN.

[2] In the IEEE format, the "implicit 1" is a 1 at the beginning of the mantissa that is not explicitly stored as part of the 32-bit representation.

## 2.2 RPC Notation

Because RPC involves computations on portions of a floating point number's representation, it is helpful to introduce some notation for describing these quantities. Assume that the RPC checker uses a fraction with *k* bits, where $k \leq 23$. Thus, an RPC checker operates on numbers that have 9+*k* bits, where the 9 corresponds to 1 sign bit and 8 exponent bits.

We use *A* and *B* to refer to the FPU's operands and *C* to refer to its (rounded) result. We use $\hat{C}$ to denote the "true" result that would be produced with unlimited precision. For example, for multiplication, $A \times B = \hat{C}$. We use the notation $A \times B \rightarrow C$ to denote that *C* is the result produced by the hardware. We denote the rounding error associated with the mantissa of $C$ ($M_C$) as $\boldsymbol{\delta_C}$. Thus, $A \times B = \hat{C} = S_C \times 2^{E_C} \times (M_C + \delta_C)$.

Assume the full-precision FPU computes $A \times B \rightarrow C$, where *A*, *B*, and *C* are full 32-bit numbers. For purposes of comparing its result, *C*, to the result of the checker, we consider the most significant 9+*k* bits of *C* and denote them as $C^H$.

The checker takes as inputs the 9+*k* most significant bits of *A* and *B*, which we denote as $A^H$ and $B^H$, respectively. It produces a 9+*k*-bit result that we denote as *C';* that is, $A^H \times B^H \rightarrow C'$. To detect errors, RPC compares *C'* to $C^H$.

The checker's inputs, $A^H$ and $B^H$, are the same as *A* and *B*, except their fractions are truncated. We denote the floating point value of $A^H$'s mantissa, including the implicit one, as $\boldsymbol{M_A^H}$. To denote the floating point value of the mantissa that gets truncated off, we use $\boldsymbol{M_A^L}$, which has value of $M_A - M_A^H$. Because $M_A$ and $M_A^H$ both contain the implicit one, their difference $M_A^L$ does not. For instance, if $M_A = 1.\{1\}_{23}$ and $k = 7$, $M_A^H = 1.\{1\}_7$ and $M_A^L = 0.\{0\}_7\{1\}_{16}$. The subscript indicates number of repetition of the bit within braces. In this case, $0.\{0\}_7\{1\}_{16}$ means 7 zeros and 16 ones following the decimal point.

The checker's output, *C'*, has a fraction and mantissa that we denote as $\boldsymbol{f'_C}$ and $\boldsymbol{M'_C}$, respectively.

# 3 RPC AND ITS CHALLENGES

At a high level, RPC is similar to dual modular redundancy. Assume the FPU performs a full-precision floating point operation on operands *A* and *B*, say $A \times B$, and obtains the full-precision result *C*. To check this result, the checker performs the same floating point operation on the inputs $A^H$ and $B^H$, where $A^H$ and $B^H$ are *9+k* most significant bits of operands *A* and *B*, respectively. The checking logic compares the *9+k*-bit result from the checker, $C' \leftarrow A^H \times B^H$, to the *9+k* most significant bits of *C* from the full-precision FPU, denoted $C^H$; if *C' significantly* differs from $C^H$, an error is reported. We illustrate this simplified version of RPC in Figure 2.

Although RPC is quite simple at this high level of abstraction, it is challenging to efficiently compare $C^H$ to *C'*. With identical dual modular redundancy, comparison between the two results can be efficiently implemented with bitwise comparison. In contrast, we need to compute the difference between the two results, because $C^H$ and *C'* can differ even in the absence of errors. However, we do not want to perform floating point subtraction, because of the time and energy required to do so.

Instead of comparing the floating point values $C^H$ and *C'* by performing a floating point subtraction, RPC compares them using **integer** subtraction. First, RPC checks whether the first bits of $C^H$ (*C*[31]) and *C'* (*C'*[31]) are equal; these bits are the signs of *C* and *C'* and they should be equal in the absence of errors. Next, RPC treats the remaining bits of $C^H$ and *C'* as *unsigned integers* and computes their integer difference (*C*[30:(23-k)]-*C'*[(7+k):0]), which we refer to as ***Diff***. The comparison logic checks whether *Diff* is within a range—bordered by an upper bound (UB) and a lower bound (LB)—that is possible in the absence of errors.

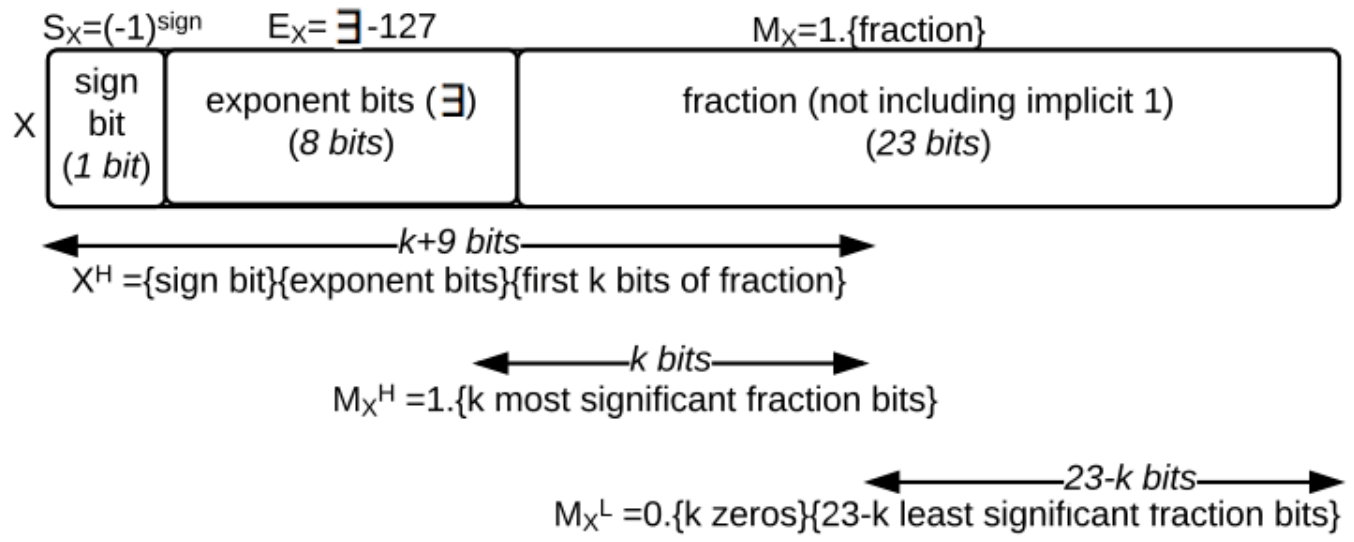

Figure 1. Representing a floating point number X

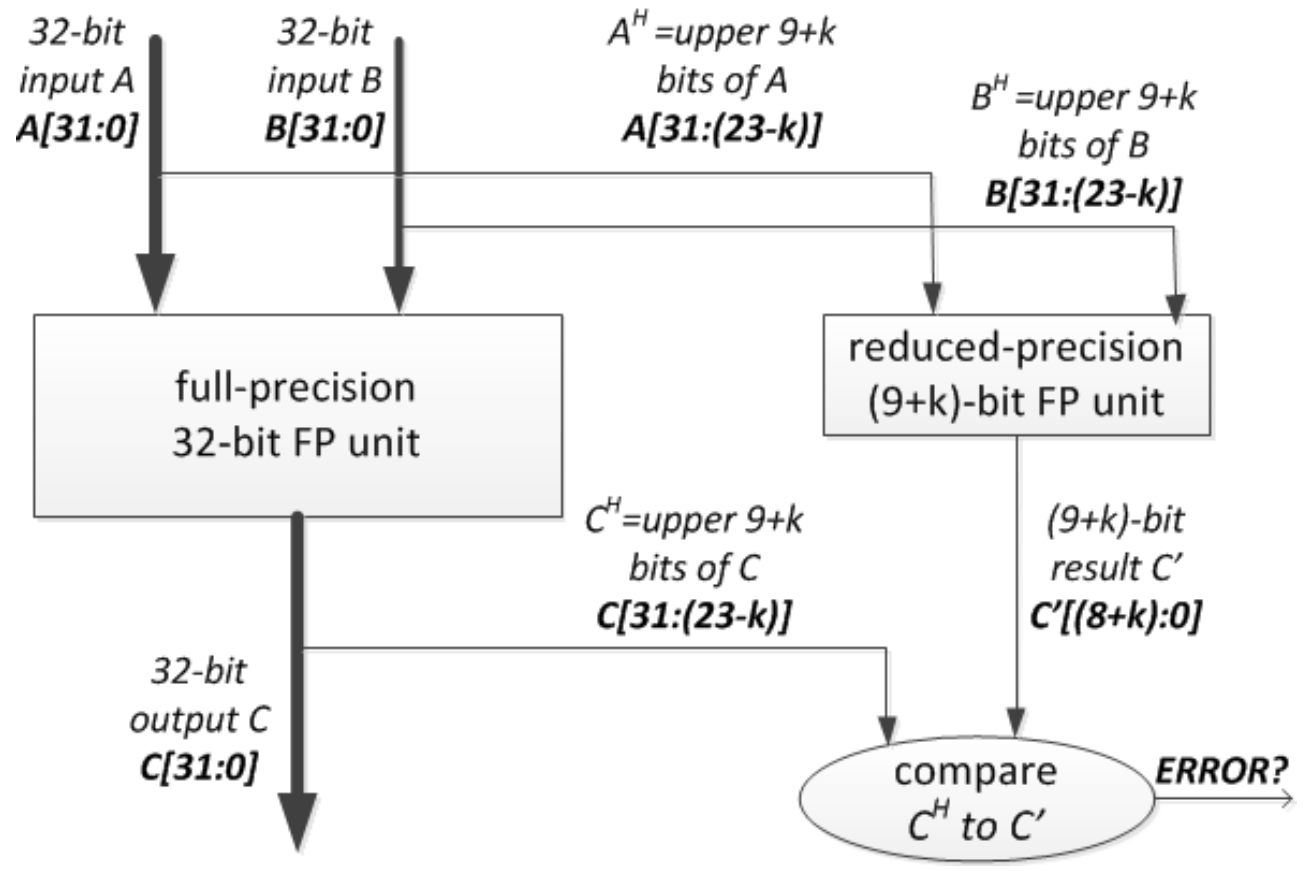

Figure 2. RPC: Big Picture

$$Diff = C[30:(23-\mathrm{k})] - C'[(7+\mathrm{k}):0]$$

$$\mathrm{NoError} = (C[31] = C'[31]) \wedge (LB < Diff < UB)$$

In Sections 5-7, we prove that, in the absence of errors:

*The bounds on Diff for addition and subtraction are [-1,1] and the [-1,3] for multiplication, division, and square root.*

If *Diff* is outside of its allowable range, then either an error has occurred in the full-precision FPU or in the checker. (The latter scenario is a false positive.) If *Diff* is non-zero but within its allowable range, then an error may or may not have occurred; if an error did occur, RPC would not detect it.

Later, we prove that the integer values of the bounds on *Diff* do not vary with *k*. However, the *significance* of these integer bounds does vary. The least significant bit of *Diff* is in the least significant bit position of the checker, which has the floating point value $2^{-k}$. Intuitively, a larger value of *k* corresponds to a smaller significance of *Diff* (i.e., $2^{-k}$ is smaller) and thus a smaller magnitude of error that will not be detected by RPC. As a result, RPC provides a tune-able tradeoff between error coverage and cost.

We devote much of this paper to proving how large the integer difference of *C'* and $C^H$ may be in the absence of errors; any difference larger than this allowable discrepancy will be detected as an error, and any difference smaller will be ignored.

Although integer subtraction is cheap and easy to implement, it introduces the following challenges in bounding the maximum allowable difference in the absence of errors.

### 3.1 Challenge #1: Unequal "Matches" of Exponents

It might appear at first that $C^H$ and *C'* can differ only in their fractions. However, they can also differ in their exponents in rare situations that must be considered. When we truncate the inputs to the checker (i.e., convert *A* and *B* into $A^H$ and $B^H$ in multiplication, for example), it leads to the absolute value of *C'* being less than or equal to the absolute value of *C*. Thus, the exponent of *C'* ($E'_C$) can be less than the exponent of $C^H$ ($E_C$). The mis-match in exponents further causes mis-alignment between the fractions of *C'* and $C^H$ for purposes of computing integer difference. This rare scenario leads to complexity in bounding *Diff*.

### 3.2 Challenge #2: The Need for "Reverse" Checking

The second challenge with RPC is the existence of certain scenarios in which direct checking is ineffective. We illustrate this challenge with a base-10 example: consider the case of $3.9999\text{x}10^{23}$-$3.9996\text{x}10^{23}$, and assume that the reduced precision checker truncates the last two digits. The (error-free) full-precision result is *C*=$3.0000\text{x}10^{19}$, and the RPC result *C'* is $3.99\text{x}10^{23}$-$3.99\text{x}10^{23}$=0. The discrepancy between $C^H$ ($3.00\text{x}10^{19}$) and *C'* in scenarios like this can be almost arbitrarily large because the result of the reduced-precision subtraction is zero. More generally, the problem exists for subtraction when both operands have the same sign and for addition when both operands have different signs. Prior work detected these scenarios and disabled checking when they occurred [1].

To overcome this challenge, we propose *reverse checking* for same-sign subtraction (SSSUB) and different-sign addition (DSADD). In the case of SSSUB, forward checking would compute $A^H - B^H \to C'$ and compare *C'* to $C^H$. Instead, with reverse checking, RPC uses the truncated result of the full-precision FPU, $C^H$, as one of it inputs and then "reverses" the computation to check it.

The operation of reverse checking depends on the signs of *A*, *B*, and *C*. There are two cases. First, if $S_A = S_B = S_C$ then the checker computes $C_H + B_H \to A'$ and then compares *A'* to $A^H$. This check avoids the problem of SSSUB because the checker is performing same-sign addition (SSADD) when computing $C_H + B_H \to A'$. In the second case, $S_A = S_B \neq S_C$, and the checker computes $A^H - C^H \to B'$ and compares *B'* to $B^H$. Because $S_C = -S_A$, the checker is performing different-sign subtraction (DSSUB) in computing $A^H - C^H \to B'$.

For division and square root, we also adopt reverse checking, but for a different reason. Both division and square root can be checked, using reverse checking, with the same reduced-precision multiplier we use to check multiplication. Reusing the multiplier checker is efficient, not least because multiplication is cheaper than division and square root. To check $A \div B \to C$, the checker performs $C^H \times B^H \to A'$ and compares $A'$ to $A^H$. To check $sqrt(B) \to C$, the checker performs $C^H \times C^H \to B'$ and compares $B'$ to $B^H$.

### 3.3 Challenge #3: the Rounding Effect

Another challenge in bounding *Diff* is that both the full-precision FPU and the reduced-precision checker produce inexact, rounded results. We must consider this rounding in two scenarios.

First, rounding in the full-precision unit might produce a carry from low-order bits that could propagate to high-order bits, leading to a mis-match between the high-order bits of $C^H$ and *C'*.

Second, rounding occurs twice in reverse checking, because one input operand used by the checker is rounded (and the result is rounded). For instance, in the case of $A \div B \to C$, even with full-precision $C \times B \neq A$. With reverse checking, we compute $C^H \times B^H \to A'$, during which we first round the input $C^H$ (by truncating *C*) and then round the result to get *A'*.

## 4 AXIOMS

In the following sections, we frequently use several axioms that derive from the format of IEEE floating point numbers.

The smallest possible value of $M_X$ is $1.\{0\}_{23} = 1$. The largest possible value of $M_X$ is $1.\{1\}_{23} = 10.\{0\}_{23} - 0.\{0\}_{22}\{1\} = 2 - 2^{-23}$.

**Axiom 1:** $1 \leq M_X \leq 2 - 2^{-23}$

Similarly, the smallest possible value of $M_X^H$ is also 1 and its largest possible value is $1.\{1\}_k = 10.\{0\}_k - 0.\{0\}_{k-1}\{1\} = 2 - 2^{-k}$.

**Axiom 2:** $1 \leq M_X^H \leq 2 - 2^{-k}$

For $M_X^L$, the smallest possible value is $0.\{0\}_{23}$ and the largest possible value is $0.\{0\}_k\{1\}_{23-k} = 0.\{0\}_{k-1}\{1\}\{0\}_{23} - 0.\{0\}_{22}\{1\} = 2^{-k} - 2^{-23}$.

**Axiom 3:** $0 \leq M_X^L \leq 2^{-k} - 2^{-23}$

For all floating point operations, the result gets rounded even in the full-precision unit. Take multiplication $A \times B \to C$ for example:

$$A \times B = \hat{C} = S_C \times 2^{E_C} \times (M_C + \delta_C) \qquad \text{Eqn. (1)}$$

Rounding error $\delta_C$ is bounded by machine epsilon ($\epsilon$). In the following derivations, we assume the default IEEE rounding mode round-to-nearest (to even on tie), in which case $\epsilon$ equals $2^{-24}$. For the other three modes (toward 0, $+\infty$, and $-\infty$), $\epsilon$ equals $2^{-23}$.

**Axiom 4:** $-2^{-24} \leq \delta_C \leq 2^{-24}$

For the reduced precision checker, the rounding error ($\delta'_C$) is still 0.5 ulp of the checker. Because the least significant bit in the checker's mantissa is $2^{-k}$, $\delta'_C$ is bounded by

**Axiom 5:** $-2^{-k-1} \leq \delta'_C \leq 2^{-k-1}$

Note that **Axiom 4** and **Axiom 5** apply to all operation types, not just multiplication, because we considered only the relationship between the result produced by a FPU and the true result with unlimited precision. Regardless of the type of operation, rounding is performed in the same manner.

An important part of bounding *Diff* is to connect the integer subtraction used to calculate *Diff* to the floating point values of its two operands. Consider the integer subtraction for computing *Diff* in forward checking.

$$Diff = C[30:(23-k)] \xrightarrow{int} C'[(7+k):0]$$

$$= \{\exists_C\}\{f_C^H\} \xrightarrow{int} \{\exists'_C\}\{f'_C\}$$

When $E_C = E'_C \rightarrow \exists_C = \exists'_C$ (common case scenario):

$$Diff = f_C^H \xrightarrow{int} f'_C$$

Because $1.\{f_C^H\} = M_C^H$ and $1.\{f'_C\} = M'_C$, the difference above equals the floating point difference between $M_C^H$ and $M'_C$ left shifted by *k* (i.e., converted to an integer).

**Axiom 6**: When $E_C = E'_C$, $Diff = \left(M_C^H \xrightarrow{float} M'_C\right) \times 2^k$

**Axiom 6** applies to reverse checking similarly, except that, instead of $C^H - C'$, *A'* or *B'* is subtracted from $A^H$ or $B^H$, respectively.

# 5 RPC FOR MULTIPLICATION

The full-precision FPU computes $A \times B \rightarrow C$:

$$(S_A \times 2^{E_A} \times M_A) \times (S_B \times 2^{E_B} \times M_B) = S_C \times 2^{E_C} \times (M_C + \delta_C)$$

During error-free execution, $S_C = S_A \times S_B$. So the signs on both sides of the equation cancel and leave:

$$M_A \times M_B \times 2^{E_A+E_B} = (M_C + \delta_C) \times 2^{E_C} \quad \text{Eqn. (2)}$$

$$M_C^H = M_A \times M_B \times 2^{E_A+E_B-E_C} - \delta_C - M_C^L \quad \text{Eqn. (3)}$$

The RPC unit computes $A^H \times B^H \rightarrow C'$:

$$(S_A \times 2^{E_A} \times M_A^H) \times (S_B \times 2^{E_B} \times M_B^H) = S'_C \times 2^{E'_C} \times (M'_C + \delta'_C)$$

Because $S_A \times S_B = S'_C$,

$$M_A^H \times M_B^H \times 2^{E_A+E_B} = (M'_C + \delta'_C) \times 2^{E'_C} \quad \text{Eqn. (4)}$$

$$M'_C = M_A^H \times M_B^H \times 2^{E_A+E_B-E'_C} - \delta'_C \quad \text{Eqn. (5)}$$

After the hardware adds *E<sub>A</sub>* and *E<sub>B</sub>* and multiplies *M<sub>A</sub>* and *M<sub>B</sub>* with integer operations, it must normalize the result. To normalize, the hardware shifts $M_A \times M_B$ until only one non-zero bit is to the left of the binary point. The number of right shifts is then added to $E_A + E_B$ to get $E_C$. By **Axiom 1**, both $M_A$ and $M_B$ are less than 2. Therefore $M_A \times M_B$ is strictly less than 4 (100 in binary), and there can be at most 2 non-zero bits to the left of the binary point. Thus, $M_A \times M_B$ is right shifted by at most one bit.

*In multiplication, the exponent of the product is equal to or is one larger than the sum of the exponents of the operands.*

We apply this rule to both the FPU and the checker:

FPU: $E_A + E_B = E_C$ or $E_A + E_B + 1 = E_C$

Checker: $E_A + E_B = E'_C$ or $E_A + E_B + 1 = E'_C$

Because operands used by the checker are truncated from operands used by the full-precision FPU, the absolute value of *C'* is less than or equal to the absolute value of *C*. So $E_C'$ must be less than or equal to $E_C$.

$$E'_C \le E_C$$

Combining the conditions above, there are three cases we need to consider.

**Common Case 1:** $E_A + E_B = E'_C = E_C$ Eqn. (6)

**Common Case 2:** $E_A + E_B + 1 = E'_C = E_C$ Eqn. (7)

**Corner Case:** $E_A + E_B + 1 = E'_C + 1 = E_C$ Eqn. (8)

The common cases are discussed in Sections 5.1 and 5.2. For both common cases, because $E_C = E'_C$, **Axiom 6** is valid, and thus we can subtract Eqn. (5) from Eqn. (3) to get:

$$Diff \div 2^k = M_C^H - M'_C$$

$$= (M_A^H M_B^L + M_A^L M_B^H + M_A^L M_B^L) \times 2^{E_A+E_B-E_C} - M_C^L + (\delta'_C - \delta_C)$$

Let $M^* = M_A^H M_B^L + M_A^L M_B^H + M_A^L M_B^L$ Eqn. (9)

$$M_C^H - M'_C = M^* \times 2^{E_A+E_B-E_C} - M_C^L + (\delta'_C - \delta_C) \quad \text{Eqn. (10)}$$

Eqn. (8) is the corner case when $E'_C \ne E_C$ so *Diff* no longer equals $(M_C^H - M'_C) \times 2^k$. We present this corner case in Appendix B.

## 5.1 Common Case 1

Using Eqn. (6), we simplify Eqn. (10) to

$$M_C^H - M'_C = \boxed{M^*}_1 + \boxed{-M_C^L}_2 + \boxed{\delta'_C - \delta_C}_3 \quad \text{Eqn. (11)}$$

Using **Axiom 2** and **Axiom 3**, we can bound boxed term 1 in Eqn. (11), which we denote <1> hereafter:

$$0 \le M^* \le (2 - 2^{-k})(2^{-k} - 2^{-23}) \times 2 + (2^{-k} - 2^{-23})^2$$

$$= (4 - 2^{-k} - 2^{-23})(2^{-k} - 2^{-23})$$

$$= 4 \cdot (2^{-k} - 2^{-23}) - 2^{-2k} + 2^{-46}$$

$$\le 4 \cdot (2^{-k} - 2^{-23})$$

Using **Axiom 3**, we bound <2>:

$$-2^{-k} + 2^{-23} < -M_C^L \le 0$$

Using **Axiom 4** and **Axiom 5**, we bound <3>:

$$-2^{-24} - 2^{-k-1} \le \delta'_C - \delta_C \le 2^{-24} + 2^{-k-1}$$

We calculate the upper bound of *Diff* by summing the upper bounds of the three terms and multiplying by $2^k$.

$$Diff = (M_C^H - M'_C) \times 2^k$$

$$\le [4 \cdot (2^{-k} - 2^{-23}) + 2^{-24} + 2^{-k-1}] \times 2^k < 4.5$$

Similarly for the lower bound on *Diff*, we have:

$$Diff = (M_C^H - M'_C) \times 2^k$$

$$\ge [-2^{-k} + 2^{-23} - 2^{-24} - 2^{-k-1}] \times 2^k > -1.5$$

The above analysis for the upper bound is simplified for purposes of providing intuition. With more sophisticated analysis (Appendix C), we can prove that the upper bound of *Diff* is strictly less than 4. The reason is because the three bounded terms above are not independent from each other and cannot reach their maximums at the same time. As *Diff* is an integer, the allowable range of values for *Diff* in this case is [-1,3].

## 5.2 Common Case 2

Using Eqn. (7), Eqn. (10) can be rewritten as:

$$M_C^H - M'_C = \boxed{M^* \times 2^{-1}}_1 + \boxed{-M_C^L}_2 + \boxed{\delta'_C - \delta_C}_3$$

Terms <2> and <3> are the same as in Section 5.1 and thus have the same bounds. Term <1> has an additional $2^{-1}$ factor not present in Section 5.1. Thus <1> is bounded by

$$0 \le M^* \times 2^{-1} \le 2 \cdot (2^{-k} - 2^{-23})$$

The upper bound on *Diff* is:

$$Diff = (M_C^H - M'_C) \times 2^k$$

$$\le [2 \cdot (2^{-k} - 2^{-23}) + 2^{-24} + 2^{-k-1}] \times 2^k \le 2.5$$

The lower bound is the same as in Section 5.1 (i.e., -1). The allowable range for *Diff* in this case is [-1,2].

> **Summary:** In this section, Appendix B, and Appendix C, we have proved that the allowable range of *Diff* for multiplication is [-1, 3].

Prior work [9] proved a looser bound of [-7,7].

# 6 RPC FOR DIVISION AND SQUARE ROOT

RPC for division and square root uses reverse checking with multiplication. We use the same reduced-precision multiplier that checks multiplication.

## 6.1 Division

To check $A \div B \rightarrow C$, the checker performs $C^H \times B^H \rightarrow A'$. In this section, we show that reverse checking of division has the same bounds on *Diff* as multiplication.

The full precision FPU computes $A \div B \rightarrow C$:

$$(S_A \times 2^{E_A} \times M_A) \div (S_B \times 2^{E_B} \times M_B) = S_C \times 2^{E_C} \times (M_C + \delta_C)$$

Because $S_A \div S_B = S_C$,

$$M_A \div M_B \times 2^{E_A - E_B} = (M_C + \delta_C) \times 2^{E_C}$$

$$M_A^H = (M_C + \delta_C) \times M_B \times 2^{E_C + E_B - E_A} - M_A^L \qquad \text{Eqn. (12)}$$

The checker computes $C^H \times B^H \rightarrow A'$:

$$(S_C \times 2^{E_C} \times M_C^H)(S_B \times 2^{E_B} \times M_B^H) = S'_A \times 2^{E'_A} \times (M'_A + \delta'_A)$$

Because $S_C \times S_B = S'_A$,

$$M_C^H \times M_B^H \times 2^{E_B + E_C} = (M'_A + \delta'_A) \times 2^{E'_A}$$

$$M'_A = M_C^H M_B^H \times 2^{E_B + E_C - E'_A} - \delta'_A \qquad \text{Eqn. (13)}$$

Using similar analysis as for multiplication, because $1 \leq M_A, M_B < 2$ according to **Axiom 1**, the FPU's quotient $M_A \div M_B$ is in the range of (0.5, 2), open interval. Therefore, the quotient needs to be left-shifted by at most 1 to be normalized such that one bit is to the left of the binary point.

*In division, the exponent of the quotient is equal to or is one less than the difference between the exponents of the dividend and divisor.*

Applying this rule to the FPU and retaining the same rule for the checker that we had for multiplication, we have:

FPU: $E_A - E_B = E_C$ or $E_A - E_B - 1 = E_C$

Checker: $E_C + E_B = E'_A$ or $E_C + E_B + 1 = E'_A$

The absolute value of *A'* is less than or equal to the absolute value of *A* due to truncating the mantissas of $M_C$ and $M_B$, and thus $E'_A \leq E_A$.

Combining the conditions above, there are again three cases we must consider:

**Common Case 1**: $E_B + E_C = E'_A = E_A$ Eqn. (14)

**Common Case 2:** $E_B + E_C + 1 = E'_A = E_A$ Eqn. (15)

**Corner Case:** $E_B + E_C + 1 = E'_A + 1 = E_A$ Eqn. (16)

The first two cases are presented in Sections 6.1.1 and 6.1.2. Because $E'_A = E_A$, $Diff = (M_A^H - M'_A) \times 2^k$ according to **Axiom 6**. The difference between Eqn. (12) and Eqn. (13) can be simplified to:

$$M_A^H - M'_A = (M_B^H M_C^L + M_B^L M_C^H + M_B^L M_C^L) \times 2^{E_C + E_B - E_A} - M_A^L$$
$$+ (\delta_C M_B \times 2^{E_C + E_B - E_A} + \delta'_A)$$

Let $M^* = M_B^H M_C^L + M_B^L M_C^H + M_B^L M_C^L$,

$$M_A^H - M'_A = M^* \times 2^{E_C + E_B - E_A} - M_A^L$$
$$+ \delta_C M_B \times 2^{E_C + E_B - E_A} + \delta'_A \qquad \text{Eqn. (17)}$$

The corner case when $E_A \neq E'_A$ can be analyzed similarly to the corner case of multiplication in Appendix B.

### 6.1.1 Common Case 1

Using Eqn. (14), Eqn. (17) can be simplified to:

$$M_A^H - M'_A = \boxed{M^*}_1 + \boxed{-M_A^L}_2 + \boxed{M_B \delta_C + \delta'_A}_3 \qquad \text{Eqn. (18)}$$

<1> and <2> have the same boundary conditions as <1> and <2> in Section 5.1, respectively. The only difference between Eqn. (18) and Eqn. (11) is that in <3> $\delta'_C - \delta_C$ becomes $M_B \delta_C + \delta'_A$. Given $M_B, \delta_C,$ and $\delta'_A$ are in the range $[1,2)$, $[-2^{-24}, 2^{-24}]$, and $[-2^{-k-1}, 2^{-k-1}]$, according to **Axiom 1**, **Axiom 4**, and **Axiom 5**, respectively, <3> now has the following bounds:

$$-2^{-k-1} - 2^{-23} < M_B \delta_C + \delta'_A < 2^{-k-1} + 2^{-23}$$

which slightly differs from the bounds on <3> in Section 5.1. This slight difference in this term does not affect the overall bounds on $M_A^H - M'_A$ as compared to Section 5.1.

$$Diff = (M_A^H - M'_A) \times 2^k$$
$$\leq [4 \cdot (2^{-k} - 2^{-23}) + 2^{-23} + 2^{-k-1}] \times 2^k < 4.5$$
$$Diff = (M_A^H - M'_A) \times 2^k$$
$$\geq [-2^{-k-1} - 2^{-23} - 2^{-k} + 2^{-23}] \times 2^k = -1.5$$

Again, the upper bound can never reach 4 for the same reason as in multiplication. Therefore, the allowable range of *Diff* is still [-1,3].

### 6.1.2 Common Case 2

Using Eqn. (15), Eqn. (17) can be simplified to:

$$M_A^H - M'_A = \boxed{M^* \times 2^{-1}}_1 + \boxed{-M_A^L}_2 + \boxed{\delta'_A + M_B \delta_C \times 2^{-1}}_3$$

The boundary condition of <3> becomes

$$-2^{-k} - 2^{-24} \leq \left(\delta'_A + M_B \delta_C \times 2^{-1}\right) \leq 2^{-k} + 2^{-24}$$

<1>, <2>, and <3> have the same boundary conditions as <1>, <2>, and <3> in Section 5.2, respectively. Therefore, both the upper and lower bounds of *Diff* are identical to those in Section 5.2.

> **Summary:** The corner case of division can be proved in the same manner as multiplication presented in Appendix B. Overall, *Diff* has the same bounded range for division as for multiplication, which is [-1,3].

## 6.2 Square Root

Checking square root is almost identical to checking division. To check $sqrt(B) \rightarrow C$, the checker performs $C^H \times C^H \rightarrow B'$. The analysis for square root is identical to the analysis for division except the two operands of the checker are the same. Therefore, the allowable range of *Diff* for square root is also [-1,3].

# 7 RPC FOR ADDITION/SUBTRACTION

Because of the cancellation problem we discussed in Section 3.2 for addition and subtraction [1], we must separately consider the following two scenarios:

- addition of same-sign operands (SSADD) or subtraction of different-sign operands (DSSUB) (Section 7.1)

- addition of different-sign operands (DSADD) or subtraction of same-sign operands (SSSUB) (Section 7.2)

## 7.1 Same-Sign Addition (SSADD) or Different-Sign Subtraction (DSSUB)

For both addition of same-sign operands ($S_A = S_B$) and subtraction of different-sign operands ($S_A \neq S_B$), RPC checks the result of the baseline FPU using forward checking, i.e., $A + B \to C$ is checked with $A^H + B^H \to C'$, and $A - B \to C$ is checked with $A^H - B^H \to C'$. These two situations can be easily converted to one another to get $A + (-B) \to C$. The resulting operation is an addition with same-sign operands (*A* and *-B*). As a result, SSADD and DSSUB have identical upper and lower bounds for *Diff*. Without loss of generality, we explain only SSADD in this section but the bounds apply equally to DSSUB.

The full-precision FPU computes $A + B \to C$:

$$(S_A \times 2^{E_A} \times M_A) + (S_B \times 2^{E_B} \times M_B) = S_C \times 2^{E_C} \times (M_C + \delta_C)$$

Because $S_A = S_B = S_C$, we have:

$$2^{E_A} \times M_A + 2^{E_B} \times M_B = 2^{E_C} \times (M_C + \delta_C) \qquad \text{Eqn. (19)}$$

$$M_C^H = M_A \times 2^{E_A - E_C} + M_B \times 2^{E_B - E_C} - M_C^L - \delta_C \qquad \text{Eqn. (20)}$$

The RPC checker computes $A^H + B^H \to C'$:

$$(S_A \times 2^{E_A} \times M_A^H) + (S_B \times 2^{E_B} \times M_B^H) = S_C' \times 2^{E_C'} \times \left(M_C' + \delta_C'\right)$$

Because $S_A = S_B = S_C'$,

$$2^{E_A} \times M_A^H + 2^{E_B} \times M_B^H = 2^{E_C'} \times (M_C' + \delta_C') \qquad \text{Eqn. (21)}$$

$$M_C' = M_A^H \times 2^{E_A - E_C'} + M_B^H \times 2^{E_B - E_C'} - \delta_C' \qquad \text{Eqn. (22)}$$

As with multiplication, division, and square root, a key aspect of our derivations is to understand the relationship between the exponents of the operands and results. To sum two numbers with the same sign, the mantissa of the operand with the smaller exponent ($E_S$) must right shift until its exponent matches the operand with the larger exponent ($E_L = \max(E_A, E_B)$). If exactly one bit is to the left of the binary point after summing the two mantissas, then the exponent of the result equals $E_L$. Else, if more than one bit is to the left of the binary point, then the sum of the mantissas must right-shift accordingly, and the exponent of the result equals $E_L$ plus the number of right shifts. The sum of the two mantissas reaches its maximum when both operands have the same exponent (i.e., no shift in mantissa of either operand), and this maximum sum of the mantissas is strictly less than 4 (100 in binary) by **Axiom 1**. As a result, the sum of the two mantissas has at most two bits to the left of the binary point, and thus the number of right shifts is either zero or one.

Based on this analysis, the exponent of the sum, $E_C$ , equals $E_L$ or $E_L$ +1. When $E_C = E_L$ , we further know that $E_A \neq E_B$ because, if $E_A = E_B$, then we can add $M_A$ and $M_B$ without shifting. The smallest value $M_A + M_B$ can be is 2 (10 in binary) because $M_A, M_B \geq 1$ according to **Axiom 1**. Thus, the sum must right shift to normalize $M_C$, which violates $E_C = E_L$.

*For addition, the exponent of the sum is equal to or one greater than the maximum of the operands' exponents ($E_L$). Furthermore, when the exponent of the sum equals $E_L$, the two operands must have unequal exponents.*

Applying this rule to the FPU and checker, we have:

FPU: $\max(E_A, E_B) = E_C \textit{ and } E_A \neq E_B$
or $\max(E_A, E_B) + 1 = E_C$

Checker: $\max(E_A, E_B) = E_C' \textit{ and } E_A \neq E_B$
or $\max(E_A, E_B) + 1 = E_C'$

The absolute value of *C'* must be less than or equal to the absolute value of *C*. Therefore,

$$E_C' \leq E_C$$

Similar to multiplication and division, there are 3 types of relationships between $E_A, E_B$, and $E_C$.

**Common Case 1:** $\max(E_A, E_B) = E_C' = E_C$ Eqn. (23)

**Common Case 2:** $\max(E_A, E_B) + 1 = E_C' = E_C$ Eqn. (24)

**Corner Case:** $\max(E_A, E_B) + 1 = E_C' + 1 = E_C$ Eqn. (25)

The common cases are considered in Sections 7.1.1 and 7.1.2. For both common cases, $E_C = E_C'$, so **Axiom 6** is valid and the difference between Eqn. (20) and Eqn. (22) can be simplified to:

$$M_C^H - M_C' = M_A^L \times 2^{E_A - E_C} + M_B^L \times 2^{E_B - E_C} - M_C^L + (\delta_C' - \delta_C) \qquad \text{Eqn. (26)}$$

We present the corner case for Eqn. (25) when $E_C \neq E_C'$ in Appendix B.

### 7.1.1 Common Case 1

In Eqn. (26), *A* and *B* are symmetric, and without loss of generality, we can explain only the case when abs(A)>abs(B) and $E_A > E_B$ because $E_A \neq E_B$.

$$E_C = E_C' = E_A > E_B \qquad \text{Eqn. (27)}$$

Eqn. (26) can be simplified to:

$$M_C^H - M_C' = \underbrace{\boxed{M_A^L + M_B^L \times 2^{E_B - E_C}}}_{1} + \underbrace{\boxed{-M_C^L}}_{2} + \underbrace{\boxed{\delta_C' - \delta_C}}_{3}$$

Because $E_C > E_B$, then $E_C \geq E_B + 1$:

$$0 < 2^{E_B - E_C} \leq 2^{-1}$$

Using **Axiom 3**, we can bound <1>:

$$0 < M_A^L + M_B^L \times 2^{E_B - E_C} \leq (2^{-k} - 2^{-23}) + (2^{-k} - 2^{-23}) \times 2^{-1} = 1.5 \times (2^{-k} - 2^{-23})$$

and we can bound <2>:

$$-2^{-k} + 2^{-23} \leq -M_C^L \leq 0$$

Using **Axiom 4** and **Axiom 5**, we can bound <3>:

$$-2^{-k-1} - 2^{-24} \leq \delta_C' - \delta_C \leq 2^{-k-1} + 2^{-24}$$

As a result, the upper bound of *Diff* is:

$$(M_C^H - M_C') \times 2^k \leq (1.5 \times (2^{-k} - 2^{-23}) + 2^{-k-1} + 2^{-24}) \times 2^k < 2$$

The lower bound of *Diff* is

$$(M_C^H - M_C') \times 2^k \geq (-2^{-k} + 2^{-23} - 2^{-k-1} - 2^{-24}) \times 2^k > -1.5$$

Therefore, the allowable range of *Diff* is [-1,1].

### 7.1.2 Common Case 2

Without loss of generality, we can again explain only $abs(A) \geq abs(B)$ and $E_A \geq E_B$.

$$E_C - 1 = E_C' - 1 = E_A \geq E_B \qquad \text{Eqn. (28)}$$

Eqn. (26) can be simplified to

$$M_C^H - M_C' = \underbrace{\boxed{M_A^L \times 2^{-1} + M_B^L \times 2^{E_B - E_C}}}_{1} + \underbrace{\boxed{-M_C^L}}_{2} + \underbrace{\boxed{\delta_C' - \delta_C}}_{3}$$

Using **Axiom 3** and Eqn. (28), we bound <1>:

$$0 < M_A^L \times 2^{-1} + M_B^L \times 2^{E_B - E_C} < (2^{-k} - 2^{-23}) \times 2^{-1} \times 2 = 2^{-k} - 2^{-23}$$

Terms <2> and <3> have the same boundary conditions as in Section 7.1.1.

So the upper bound of *Diff* is:

$$(M_C^H - M_C') \times 2^k \le (2^{-k} - 2^{-23} + 2^{-k-1} + 2^{-24}) \times 2^k < 1.5$$

The lower bound of *Diff* is still greater than -1.5. Therefore, the allowable range of *Diff* is [-1,1].

## 7.2 Different-Sign Addition (DSADD) or Same-Sign Subtraction (SSSUB)

Recall that our simple example in Section 3.2 was a SSSUB that motivated reverse checking. The situation occurs when $E_A = E_B$ and $M_A^H = M_B^H$, but $M_A^L \neq M_B^L$. In this situation, the result of $A - B \rightarrow C$ computed by the full-precision FPU has a non-zero value, whereas $A^H - B^H \rightarrow C'$ computed by the checker (i.e., with forward checking) equals zero. As a result, *Diff = C[30:(23-k)]-C'[(7+k):0]* can become arbitrarily large even in the absence of errors. DSADD is equivalent to SSSUB and thus suffers from the same problem in forward checking.

To resolve the challenge, the checker applies reverse checking using $C^H$ as one of its operands. In this way, we use SSADD/DSSUB to check SSSUB/DSADD.

For SSSUB: $A - B \rightarrow C$, $S_A = S_B$:

- If $S_A = S_B = S_C$, check $C^H + B^H \rightarrow A'$ (SSADD).
- If $S_A = S_B \neq S_C$, check $A^H - C^H \rightarrow B'$ (DSSUB).

Similarly, for DSADD $A + B \rightarrow C, S_A \neq S_B$:

- If $S_A = S_C \neq S_B$, check $C^H - B^H \rightarrow A'$(DSSUB).
- If $S_A \neq S_C = S_B$, check $C^H - A^H \rightarrow B'$(DSSUB).

In the following section, we analyze the first scenario above: SSSUB with $S_A = S_B = S_C$. The other cases can be verified in a similar manner.

The full-precision FPU computes $A - B \rightarrow C$:

$$(S_A \times 2^{E_A} \times M_A) - (S_B \times 2^{E_B} \times M_B) = S_C \times 2^{E_C} \times (M_C + \delta_C)$$

Because $S_A = S_B = S_C$ in this scenario:

$$M_A - 2^{E_B - E_A} \times M_B = 2^{E_C - E_A} \times (M_C + \delta_C)$$

$$M_A^H = M_B \times 2^{E_B - E_A} + M_C \times 2^{E_C - E_A} - M_A^L + \delta_C \times 2^{E_C - E_A}$$

The checker performs $C^H + B^H \rightarrow A'$:

$$(S_C \times 2^{E_C} \times M_C^H) + (S_B \times 2^{E_B} \times M_B^H) = S_A' \times 2^{E_{A'}} \times (M_A' + \delta_A')$$

Because $S_A' = S_B = S_C$,

$$2^{E_C - E_A'} \times M_C^H + 2^{E_B - E_A'} \times M_B^H = M_A' + \delta_A'$$

$$M_A' = M_B^H \times 2^{E_B - E_A'} + M_C^H \times 2^{E_C - E_A'} - \delta_A'$$

To understand the relationships between $E_A$, $E_B$, and $E_C$ in SSSUB (which also apply to DSADD), consider $A - B = \hat{C}$, where $\hat{C}$ is the true $C$ with unlimited precision. Then $B + \hat{C} = A$. From our previous analysis in Section 7.1, we have:

$$E_A = \max(E_B, E_{\hat{C}}) \textit{ and } E_B \neq E_{\hat{C}} \quad (1)$$

or $E_A = \max(E_B, E_{\hat{C}}) + 1$ (2)

From **Axiom 4**,

$$E_C = E_{\hat{C}} \text{ if } \hat{C} \text{ rounds down} \quad (3)$$

or $E_C = E_{\hat{C}} + 1$ if $\hat{C}$ rounds up (4)

We now look at the "cross-product" of possible scenarios: (1&3), (2&3), (1&4), (2&4):

(1&3) $E_A = \max(E_B, E_C), E_B \neq E_C$

(2&3) $E_A = \max(E_B, E_C) + 1$

The other two scenarios are either impossible or subsumed by a previous scenario.

Scenario (1&4) differs from Scenario (1&3) only when $E_{\hat{C}} > E_B$, so $E_A = \max(E_B, E_C) = E_{\hat{C}}$ in (1). Then because of (4), $E_{\hat{C}} = E_C - 1$. So $E_A = E_C - 1 = \max(E_B, E_C)$. Thus, Scenario (1&4) is impossible because it implies that *abs(C)>abs(A)*, which is impossible after *B* is subtracted from *C*, given that $S_A = S_B$.

Scenario (2&4) differs from (2&3) only when $E_{\hat{C}} > E_B$. Then $E_A = E_{\hat{C}} + 1 = (E_C - 1) + 1 = E_C$. Because $E_C > E_{\hat{C}} > E_B$, $E_A = \max(E_B, E_C)$, *so* $E_A = \max(E_B, E_C)$, which is the same as scenario (1&3).

*In subtraction, the exponent of the minuend (A) is equal to or one greater than $E_L$, which is the maximum of the exponents of the subtrahend (B) and difference (C). Furthermore, when the exponent of the minuend equals $E_L$, the subtrahend and the difference must have unequal exponents.*

Applying this rule to the FPU and applying the rule for addition to the checker, we have:

FPU: $E_A = \max(E_B, E_C)$ and $E_B \neq E_C$
or $E_A = \max(E_B, E_C) + 1$

Checker: $E'_A = \max(E_B, E_C)$ and $E_B \neq E_C$
or $E_A' = \max(E_B, E_C) + 1$

Because of truncation in operands *B* and *C*, $E_A' \le E_A$.
Therefore,

**Common Case 1**: $\max(E_B, E_C) = E_A' = E_A, E_B \neq E_C$ Eqn. (29)

**Common Case 2:** $\max(E_B, E_C) + 1 = E_A' = E_A$ Eqn. (30)

**Corner Case:** $\max(E_B, E_C) + 1 = E_A' + 1 = E_A$ Eqn. (31)

We analyze Eqn. (29) and Eqn. (30) in Sections 7.2.1 and 7.2.2. Because $E_A = E_A'$ in these common cases, **Axiom 6** is valid and the difference between $M_A^H$ and $M_A'$ is

$$M_A^H - M_A' = \underbrace{\boxed{M_B^L \times 2^{E_B - E_A} + M_C^L \times 2^{E_C - E_A}}}_{1} + \underbrace{\boxed{-M_A^L}}_{2} + \underbrace{\boxed{\delta_A' + \delta_C \times 2^{E_C - E_A}}}_{3} \quad \text{Eqn. (32)}$$

Eqn. (31) can be proved in similar way as the corner case of SSADD (Appendix A).

### 7.2.1 Common Case 1

In term <1>, because of Eqn. (29), either $2^{E_B - E_A}$ or $2^{E_C - E_A}$ must equal 1 because $\max(E_B, E_C) = E_A$. The other is less than or equal to $2^{-1}$ because $E_B \neq E_C$. As a result, <1> is bounded by:

$$0 \le M_B^L \times 2^{E_B - E_A} + M_C^L \times 2^{E_C - E_A} \le (2^{-k} - 2^{-23}) + (2^{-k} - 2^{-23}) \times 2^{-1} = 1.5 \times (2^{-k} - 2^{-23})$$

For the same reason, $2^{E_C - E_A}$ in <3> is smaller than or equal to 1. So we can bound <3>:

$$-2^{-k-1} - 2^{-24} \le \delta_C \times 2^{E_C - E_A} + \delta_A' \le 2^{-k-1} + 2^{-24}$$

Term <2> is bounded to $[-2^{-k} + 2^{-23}, 0]$ according to **Axiom 3**. As <1>, <2>, and <3> have the same boundary conditions as <1>, <2>, and <3> in Section 7.1.1, respectively, the bounds on *Diff* are identical to those in Section 7.1.1. Therefore, *Diff* is in the range [-1,1].

### 7.2.2 Common Case 2

Because of Eqn. (30), $0 < 2^{E_B - E_A}, 2^{E_C - E_A} \le 2^{-1}$ in <1>. As a result, <1> is bounded by:

$$0 \le M_B^L \times 2^{E_B - E_A} + M_C^L \times 2^{E_C - E_A}$$
$$\le (2^{-k} - 2^{-23}) \times 2^{-1} + (2^{-k} - 2^{-23}) \times 2^{-1}$$
$$= 2^{-k} - 2^{-23}$$

In addition to the boundary condition from **Axiom 4** and **Axiom 5**, <3> is bounded by

$$-2^{-k-1} - 2^{-25} < \delta_C \times 2^{E_C - E_A} + \delta_A^{'} < 2^{-k-1} + 2^{-25}$$

Boxed term <3> has smaller bounds than <3> in Section 7.1.2. Furthermore, <1> and <2> have same upper and lower bounds with <1> and <2> in Section 7.1.2, respectively. Therefore, *Diff* is in the range [-1,1].

**Summary:** Along with analysis in Appendix A, the analysis in section 7 shows the allowable range of *Diff* for addition and subtraction are [-1,1].

Prior work [1] proved a looser bound of [-2,1].

## 8 MAXIMUM PERCENTAGE ERROR

With bounds on *Diff*, we can now approximate the Maximum Percentage Error (MPE) for errors undetected by RPC. This derivation is an approximation in that it applies to all common situations but not some of the extremely rare corner cases. This section shows that the undetected errors are small. We compute the percentage error as:

$$\%\text{Error} = \left|\frac{\text{Correct Result} - \text{Erroneous Result}}{\text{Correct Result}}\right| \times 100\%$$

We have proven that *Diff* is in the range of [-1,1] for addition/subtraction and [-1,3] for multiplication, division, and square root. For the common cases, *Diff* corresponds to the difference between the fractions of the baseline FPU and the checker.

We represent an erroneous computation with an "E" above the right-arrow and an erroneous result as $C_e$.

### 8.1 Forward Checking

In forward checking, the FPU and the checker compute

$$A \text{ } op \text{ } B \xrightarrow{E} C_e \text{ and } A^H \text{ } op \text{ } B^H \to C'$$

$$\text{MPE} = \left|\frac{C - C_e}{C}\right| = \left|\frac{M_C - M_C^e}{M_C}\right| \cong \left|\frac{M_C' - (M_C^e)^H}{M_C}\right| \cong \left|\frac{Diff \times 2^{-k}}{M_C}\right|$$

RPC only misses an error when *Diff* is within the derived bounds. Also $M_C \ge 1$, so

$$\text{MPE} \le \max(|Diff|) \times 2^{-k} \times 100\%$$

### 8.2 Reverse Checking

In reverse checking of *division* (which also applies to *square root*), *Diff* represents the difference in the fractions of $A^H$ and *A'*. However, what we are really interested in is the difference between $C_e$ and *C*. We know that:

$$A \div B \xrightarrow{E} C_e \text{ and } C_e^H \times B^H = A'$$

The percentage error is:

$$\left|\frac{C - C_e}{C}\right| \cong \left|\frac{A^H \div B^H - A' \div B^H}{C}\right| \cong \left|\frac{(A^H - A') \div B^H}{C}\right|$$
$$= \frac{(Diff \times 2^{-k} \div M_B) \times 2^{E_A - E_B}}{M_C \times 2^{E_C}}$$

As $E_A - E_B \ge E_C$ and $M_B M_C \ge 1$, MPE is loosely bounded by

$$MPE \le \frac{(Diff \times 2^{-k})}{M_B M_C} \le \max(|Diff|) \times 2^{-k} \times 100\%$$

In reverse checking of *SSSUB (also applies to DSADD)*,

$$A - B \xrightarrow{E} C_e \text{ and } C_e^H + B^H = A'$$

Then

$$\left|\frac{C - C_e}{C}\right| \cong \left|\frac{A - B - (A' - B^H)}{C}\right|$$

When $E_A \cong E_C \gg E_B$, $|C - C_e| \cong |A^H - A'| = |M_A^H - M_A'| \times 2^{E_C}$

$$\text{MPE} = \max\left|\frac{M_A^H - M_A{}'}{M_C}\right| \times 100\% \le \max(|Diff|) \times 2^{-k} \times 100\%$$

When $E_B \cong E_C \gg E_A$, $|C - C_e| \cong B^L = M_B^L \times 2^{E_C}$

$$\text{MPE} = \max\left|\frac{M_B^L}{M_C}\right| \times 100\% < 2^{-k} \times 100\%$$

### 8.3 Summary

Because $\max(|Diff|)$ for SSSUB is just 1, MPE is approximately $\max(|Diff|) \times 2^{-k} * 100\%$ for all operations. This approximate analysis does not completely bound the maximum percentage error uncaught by RPC, but rather shows that, under most circumstances, undetected errors have very small percentage errors. For example, the MPE for a 16-bit checker (k=7) is 0.78% for addition and subtraction, and 2.34% for the other operations.

## 9 IMPLEMENTATION AND INTEGRATION INTO PROCESSOR CORE

We implemented RPC as an extension of a floating point unit developed by Kwon et al. [3]. Our implementation of RPC includes the following four components: one *k*+9-bit adder/subtractor, one *k*+9-bit multiplier, logic to determine how to perform the checking, and buffers to hold operands and results until checking can be performed.

### 9.1 Handling Floating Point Exceptions

When certain floating point exceptions occur, RPC is unable to check for errors. The reason for this limitation is that the assumptions we make about rounding error are valid only when the FPU does not encounter the following exceptions: overflow, underflow, invalid, and divide-by-zero. In these circumstances, the result of the FPU is formatted to positive/negative infinity, zero or denorm, NaN, and positive/negative infinity, accordingly. Therefore, our checker is suppressed when these rare cases are detected.

### 9.2 Performance Impact

RPC can impact performance if the processor core is waiting for a result to be checked before it can commit a floating point instruction. For some floating point operations, we can perform checking in parallel with the full-precision FPU, but reverse checking does not permit this parallelism. If we try to begin each check as soon as it has its operands, we could encounter a situation in which a reverse check and a forward check need to use the same checker at the same time. To avoid the complexity of managing these situations, we choose to have the checker always wait until the FPU has completed (even for forward checking).

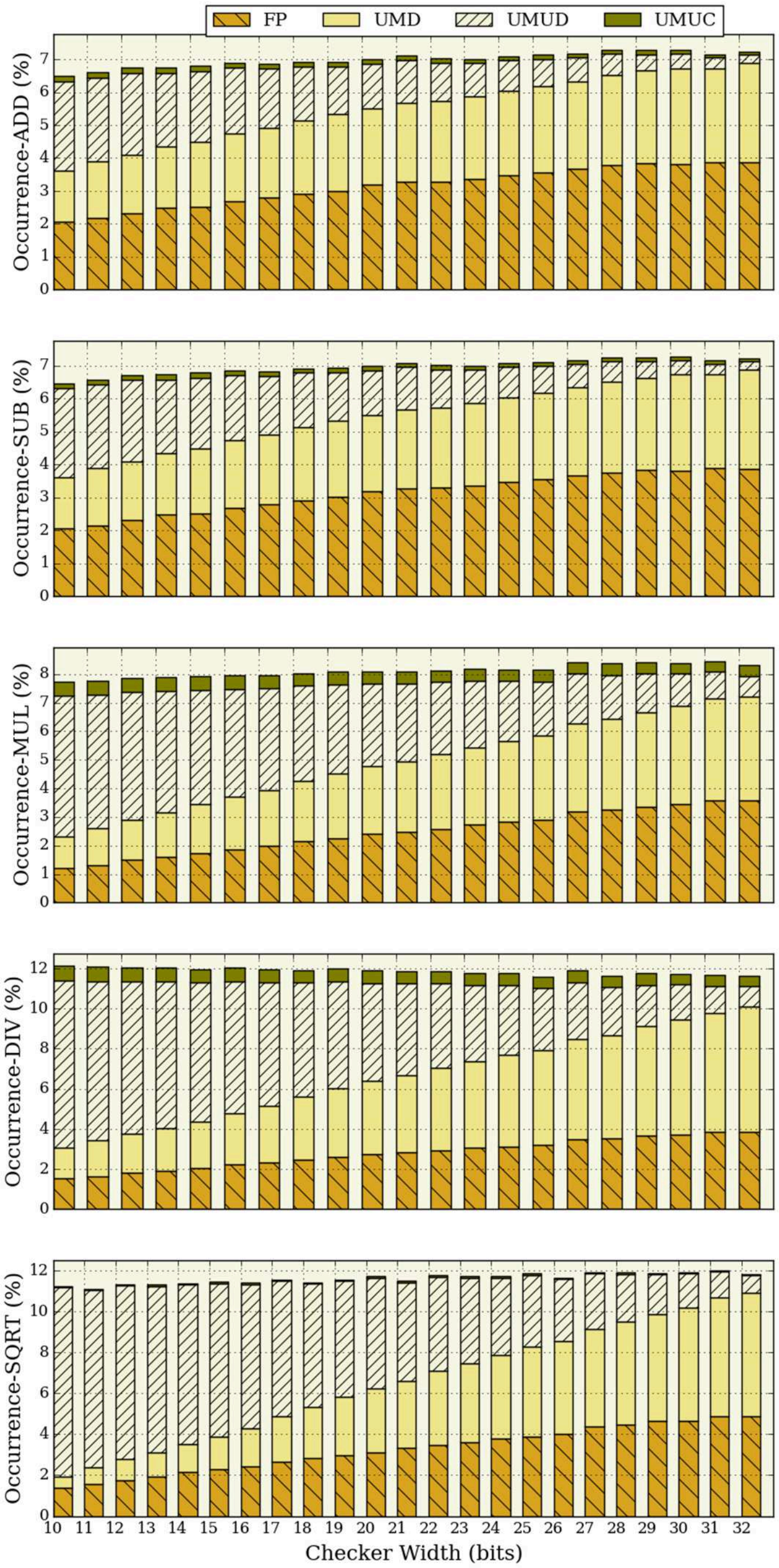

Figure 3. Classification of Injected Errors

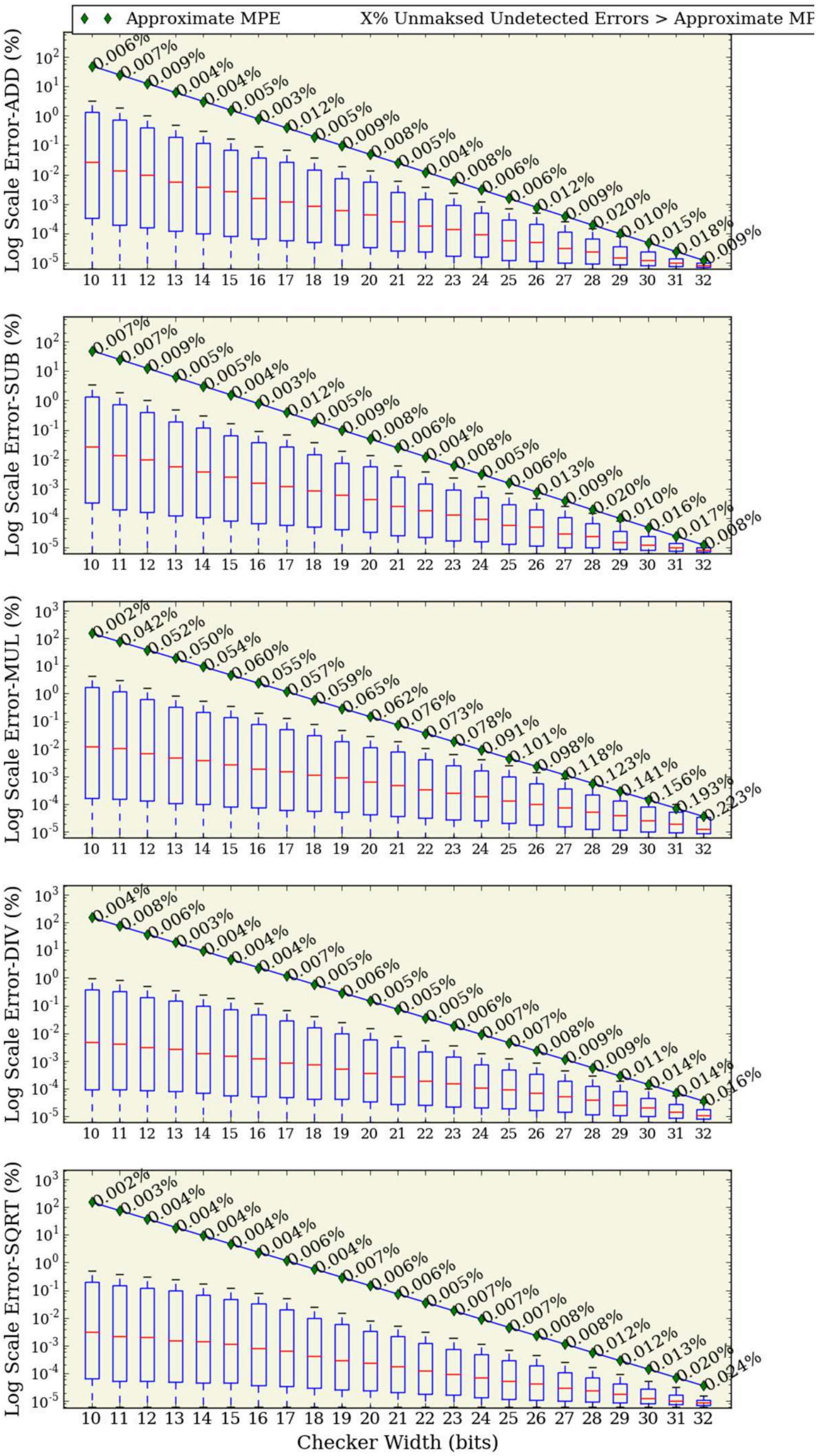

Figure 4. Distribution of Size of UMUD Errors

As the checker is modified based on the baseline FPU, it should have at least the same throughput and thus the performance impact of RPC is limited to the extra reorder buffer (ROB) pressure incurred by floating point instructions that are ready to commit but have not yet been checked. The performance penalty due to this extra ROB pressure depends on the microarchitecture and the software running on it, but we do not expect it to be large.

## 10 Error Detection Coverage

The primary goal of our experimental evaluation is to determine the error detection coverage of RPC.

### *10.1 Error Injection Methodology*

To evaluate the error detection coverage, we ran an extensive set of error injection experiments. In each experiment, we forced a single stuck-at error on a different wire in the flattened netlist that includes the full-precision FPU and all of the RPC hardware. Note that a single stuck-at error on a wire can often cause multiple errors downstream of the injected error, due to fan-out, and thus injecting errors at this low level is far more realistic than injecting errors in a microarchitectural simulator [4].

For every wire in the netlist, we ran 2000 experiments. Half of the experiments inject a stuck-at-0 error on the wire and test 1000 different inputs; the other half inject a stuck-at-1 on the wire and test the same 1000 random inputs.

We considered injecting transient errors, but a very large fraction of transient errors were masked (i.e., had no impact on execution). Masking is quite common, as in prior work in error injection [8], because each error affects the results for only a subset of all possible inputs, and often these subsets are tiny. Transient errors, in particular, are masked with very high probability. To obtain statistically significant data in a tractable amount of time, we injected permanent stuck-at errors, which are less likely to be masked. Moreover, a floating point operation is a relatively short latency event, which makes transient errors similar to permanent errors.

### *10.2 Results*

In Figure 3, we show our results by classifying errors into categories, and we present a separate graph for each of the five arithmetic operations (add, sub, mul, div, sqrt). The figures focus on the errors that are *unmasked*, i.e., those errors that affect the result of the FPU and/or the result of the checker.

Among the unmasked errors, we classify the errors into four categories:

- unmasked and detected (UMD) – our desired outcome
- unmasked undetected (UMUD) – a silent data corruption, which is the worst outcome

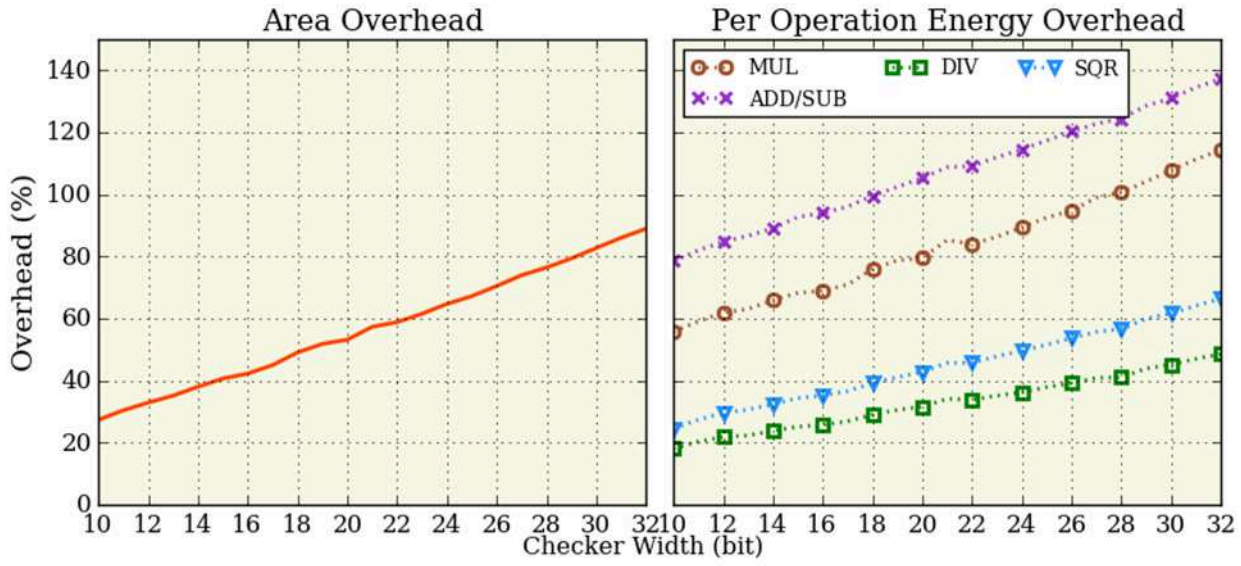

(a) (b)

Figure 5. Area and Energy Overheads

- unmasked unchecked (UMUC) –unusual corner case when RPC suppresses checking
- false positive (FP) – error affected checker and is detected even though result of FPU is correct

The first observation from Figure 3 is that there is a clear trade-off between error detection capability and cost. As the checker is made wider, there are fewer undetected errors.

A corollary to this first result is that there are also more false positives, because the checker is larger and thus more liable to be the victim of an injected error. Notice that there are fewer false positives in division and square root than addition, subtraction, and multiplication; this is because the divider and square root unit are significantly larger than the checker multiplier. Unlike a dual modular redundancy scheme (i.e., simply replicate the unit to be checked), which has 50% false positives, RPC can reduce the fraction of false positives by having a smaller checker.

From Figure 3, it appears RPC misses a large portion of unmasked errors at small checker width (such as 16, *i.e.*, 7 bits mantissa). However, Figure 4 shows that **RPC detects the vast majority of the *large* errors.**

Figure 4 presents the distribution of the *size* of the unmasked undetected (UMUD) errors. In these graphs, the top, middle, and bottom bars of the vertical box indicate the first quartile, median, and second quartile of percentage errors, respectively. The top whisker defines a 1.5 inter-quartile range away from the first quartile, and data above here are defined as outliers. The connected green diamonds mark the approximate MPE computed in Section 8. The percentage number above the diamonds is the percentage of UMUD errors that has percentage error larger than the approximate MPE. An error larger than the approximate MPE includes the corner case scenarios and faults in a handful of components that, when faulty, can cause extremely strange behavior. For example, an error in the wire that determines whether the output should be formatted in an exceptional way (e.g., as a NaN) can cause the output to differ dramatically from the correct result and cause an outlier in the percentage error calculation.

We observe that only a very small fraction of UMUDs have percentage errors that are not bounded by the MPE. Overall, the percentage errors of UMUDs are extremely small; the median is less than 0.1% even with a minimally-sized 1-bit mantissa (10-bit checker) for all operations. At checker width 16, only 0.003%~0.055% undetected errors have percentage error larger than the approximate MPE (0.78% for addition and subtraction, and 2.34% for the other operations). This means a vast majority of those undetected errors in Figure 3 at width 16 are very small.

## 11 Area, Power, and Energy Overheads

We now evaluate the area and power overheads of our RPC implementation.

### *11.1 Area*

We used Synopsys CAD tools to layout the FPU, with and without RPC, in 45nm technology [7]. The results in Figure 5 (a) show that RPC's area overhead ranges from about 30% to about 90%, depending on the checker's width. This overhead is far less than that of simple dual modular redundancy (100%).

### *11.2 Energy*

We compute the per operation energy overhead using the power and latency results from the CAD tools. As with the area analysis, we determine the dynamic and static power overheads of RPC with post-synthesis gate-level power estimation. After laying out the circuitry, we obtained the parasitic resistances and capacitances and back-annotated the circuits with them.

To determine the power, we feed the synthesized module with 1000 random inputs for each operation at all checker widths to acquire switching activities of each wire/cell, which is further used to compute both dynamic and static power. To minimize power consumed by buffers for *reverse checking*, we use clock gating to optimize our design.

Our experiments show that static power comprises roughly 20% of overall power and, unlike dynamic power, is relevantly stable across different operations and checker widths. Figure 5(b) shows the per operation energy overhead of 5 operations. Addition and subtraction suffer the most from energy overhead because they are relatively cheap operations. Hence, the checking logic incurs nontrivial overhead. However, before width 18, the overhead is still cheaper than dual modular redundancy. Multiplication incurs somewhat less energy overhead, reaching 100% at width 28. For the most expensive operations, division and square root, as their checker multiplier is significantly smaller in size, RPC is very efficient and has less than 70% overhead at full checker width. Note that the overheads in Figure 5(b) include the impact of the checker, logic, and buffer over each single operation. However, the checker adder is shared between addition and subtraction; the checker multiplier is shared among the other 3 operations; and logic and buffers are shared among all operations. So the overheads do not sum up for a mixture of 5 operations.

## 12 Related Work

The most related work is the prior work on RPC, which we have already discussed [1][9]. In addition, there have been other proposed schemes for detecting errors in floating point hardware.

Lipetz and Schwarz [5] propose residue checking, which is complete and cost-efficient, in principle, but we cannot resolve how their scheme handles the issue of rounding. We speculate that, based on an IBM patent [2], the FPU passes rounding information to the residue checker, but such a design would be unable to detect errors in rounding logic.

Maniatakos et al. [6] propose a low-cost scheme in which the hardware checks only the exponents of floating point operations. Like us, they seek to detect all large errors, but checking only exponents misses many errors that we consider too large to be acceptable.

## 13 Conclusions

In this work, we have applied RPC to an entire floating point unit. We have comprehensively analyzed its ability to detect errors and the size of errors that can go undetected. Based on our results, we believe that RPC is an attractive method for detecting errors in FPUs. We are unaware of other approaches for detecting errors in FPUs that can be tuned to trade off error detection coverage versus cost.

## Acknowledgments

This material is based on work supported by the National Science Foundation under grant CCF-111-5367.

## *Appendix A: $E_C$=$E_C$'+1 for SSADD*

In this section we consider the corner case scenario for SSADD Eqn. (25). When $E_C = E_C' + 1$, **Axiom 6** is no longer valid. To derive the upper bound of *Diff* in this case, we need to understand the difference (we denote as *D*) between what is computed by the FPU and the checker and under what circumstance *D* causes $E_C \neq E_C'$. To simplify the analysis, we introduce the notation */X/* to denote the integer represented by the concatenation of the exponent and mantissa bits of *X* without sign bits. So for example, $C[30{:}(23-k)]$ is $|C^H|$ and $C'[(7+k){:}0]$ is $|C'|$.

We know that $|C'| \leq |C|$ because of truncation in operands in computing C'. As we care only about the magnitude of the difference in this case, we define *D = /C/-/C'/*. Then */C/ = /C'/ + D.* (Note that *D* is not *Diff*, which is */C^H/-/C'/*.) *D* is the difference between the left-hand sides of Eqn. (19) and Eqn. (21):

$$D = M_A^L \times 2^{E_A} + M_B^L \times 2^{E_B}$$

Let $D = 2^{E_D} \times M_D$ and let $E_D = E_C'$. Then

$$M_D = M_A^L \times 2^{E_A - E_C'} + M_B^L \times 2^{E_B - E_C'}$$

$E_C'$ differs from $E_C$ only when $M_C' + M_D \geq 2$. Then to normalize the sum, $M_C' + M_D$ must right shift, resulting in $E_C = E_C' + 1$. From Eqn. (25), $E_C' = \max(E_A, E_B), E_A \neq E_B$. Therefore, $2^{E_A - E_C'} \leq 1$ and $2^{E_B - E_C'} \leq 1$, but they cannot both equal 1 at the same time. Thus, we have:

$$M_D < M_A^L + M_B^L < 2 \times 2^{-k}$$

Therefore, $M_D$ must be in the form of $0.\{0\}_{k-1}\{x \ldots\}$, where each *x* represents an unknown bit and can be either zero or one. Each *x* is independent from each other and they may or may not hold the same value. $M_C'$ is all zero after its first *k* bits to the right of its binary point, i.e., it is in the form $1.\{x\}_k\{0 \ldots\}$

$$\begin{array}{lc|c|c} M_D: & 0. & 000 \ldots 000x & x \ldots x \\ M_C': & 1. & xxx \ldots xxxx & 0 \ldots 0 \end{array}$$

*k* bits

If the sum of the two ever reaches 2, the first $x$ in the fraction of $M_D$ must equal 1 and all $x$ in $M'_C$ must equal 1. Their sum must be in the form of $10.\{0\}_k\{x\ldots\}$ before being normalized. After being normalized, $M_C$ must be in the form of $1.\{0\}_{k+1}\{x\ldots\}$.

$$\begin{array}{rr|c|l} M_D: & 0. & 000\ldots0001 & x\ldots x \\ +\ M'_C: & 1. & 111\ldots1111 & 0\ldots0 \\ \hline M_D+M'_C: & 10. & 000\ldots0000 & x\ldots x \\ M_C: & 1. & 000\ldots0000 & 0x\ldots x \end{array}$$

*k* bits

So when $E_C \neq E'_C$, the fractions (not including the implicit 1) must be $f_C^H = \{0\}_k, f'_C = \{1\}_k$. Now consider $Diff = |C^H| - |C'|$:

$$\frac{\begin{array}{r}|C^H|\\ -\ |C'|\end{array}}{Diff} \rightarrow \frac{\begin{array}{r}\{\exists_C\}\{f_C^H\}\\ -\ \{\exists'_C\}\{f'_C\}\end{array}}{(E_C - E'_C)\times 2^k + f_C^H - f'_C} \rightarrow \frac{\begin{array}{r}\{\exists'_C+1\}\{00\ldots0\}\\ -\ \{\exists'_C\}\{11\ldots1\}\end{array}}{1}$$

where $\exists_C = E_C + 127$ and $\exists'_C = E'_C + 127$. Therefore, for the corner case Eqn. (25), *Diff* equals 1, which is within the range [-1,1].

## Appendix B: $E_C = E_C' + 1$ for Multiplication

This corner case of multiplication, Eqn. (8), follows a similar analysis as in Appendix A. In this case, *D* is the difference between the left hand sides of Eqn. (2) and Eqn. (4).

$$D = M^* \times 2^{E_A+E_B}$$

where M* is defined in Eqn. (9). From Eqn. (8), we know that $E_A + E_B = E'_C$. Let $D = 2^{E_D} \times M_D, E_D = E'_C$. Then we have $M_D = M^*$. In Section 5.1, we bound $M^*$ in <1>. So $M_D < 4 \times 2^{-k}$ and is thus in the form of $0.\{0\}_{k-2}\{x\ldots\}$. We also know that $M'_C$ has only *k* mantissa bits, so it is of the form $1.\{x\}_k\{0\ldots\}$.

$$\begin{array}{rr|c|l} M_D: & 0. & 000\ldots00xx & x\ldots x \\ M'_C: & 1. & xxx\ldots xxxx & 0\ldots0 \end{array}$$

*k* bits

The first (*k*-2) values of *x* in the fraction of $M'_C$ must be 1 for the sum to reach 2. The sum of the last two *x* bits of $M_D$ (denoted as *ab*) and last two *y* bits of $M'_C$ (denoted as *cd*) must be greater than or equal to $4_{10} = 100_2$ to produce the carry one.

$$\begin{array}{rr|c|l} M_D: & 0. & 000\ldots00ab & x\ldots x \\ +\ M'_C: & 1. & 111\ldots11cd & 0\ldots0 \\ \hline M_D+M'_C: & 10. & 000\ldots00ef & x\ldots x \\ M_C: & 1. & 000\ldots000e & fx\ldots x \end{array}$$

*k* bits

Now *Diff* can be evaluated in the following steps:

$$\frac{\begin{array}{r}|C^H|\\ -\ |C'|\end{array}}{Diff} \rightarrow \frac{\begin{array}{r}\{\exists_C\}\{0\ldots00e\}\\ -\ \{\exists'_C\}\{1\ldots1cd\}\end{array}}{2^k - \{1\}_{k-2}\{0\}_2 + \{0e - cd\}} \rightarrow \frac{\begin{array}{r}\{\exists'_C+1\}\{0\ldots0\}\{0e\\ -\ \{\exists'_C\}\{1\ldots1\}\{cd\}\end{array}}{\{1\}\{0e\text{-}cd\}=10e\text{-}cd}$$

Consider all possible value of *ab*, *cd*, and *ef*, where *ab* + *cd* = 1*ef*. The two-bit values *ab* and *cd* are both less than or equal to $11_2 = 3_{10}$. Therefore, *ef* cannot be $11_2$, as that requires $ab + cd = 111_2 = 7_{10}$. Thus, the possible values of *ab*, *cd*, *ef*, and *Diff* are as follows:

| *ab* | *cd* | *1ef* | *10e* | *Diff* = 10*e* - *cd* |
|---|---|---|---|---|
| 01 | 11 | 100 | 100 | $1_2 = 1_{10}$ |
| 11 | 01 | | | $11_2 = 3_{10}$ |
| 10 | 10 | | | $10_2 = 2_{10}$ |
| 10 | 11 | 101 | | $1_2 = 1_{10}$ |
| 11 | 10 | | | $10_2 = 2_{10}$ |
| 11 | 11 | 110 | 101 | $10_2 = 2_{10}$ |

Therefore, *Diff* is still within the range [-1,3] of multiplication.

## Appendix C: Upper Bound of Diff in Multiplication

In Section 5.1 we proved the upper bound of *Diff* is less than 4.5, indicating that this integer difference can be as large as 4. However, after billions of simulations, we never observed *Diff* equaling 4. In this section, we introduce a more sophisticated proof that bounds *Diff* to be less than 4 (*i.e.* $M_C^H - M'_C < 4 \times 2^{-k}$) and thus this integer difference can be at most 3. This proof also implies that the upper bounds on *Diff* for division and square root are also less than 4.

To bound *Diff*, we must understand how $M'_C$ differs from $M_C^H$. Thus, we need to understand how rounding happens in the checker and how rounding in the FPU affects the most significant *k* bits in $M_C$ (i.e. $M_C^H$).

In the following analysis, we denote the unrounded (unlimited precision) result of the FPU's mantissa as $M_{\hat{C}}$ and the unrounded result of the checker's mantissa as $M'_{\hat{C}}$. With unlimited precision, $M_{\hat{C}}$ has more than 23 bits of fraction (i.e., bits after the binary point), and $M'_{\hat{C}}$ has more than *k* bits of fraction even though their input operands have only 23 and *k* bits of fractions, respectively. We split both $M_{\hat{C}}$ and $M'_{\hat{C}}$ into High and Low terms: $M_{\hat{C}} = M_{\hat{C}}^H + M_{\hat{C}}^L$ and $M'_{\hat{C}} = (M'_{\hat{C}})^H + (M'_{\hat{C}})^L$, where $M^H$ refers to the implicit 1 and the first *k* bits of the fraction; $M^L$ refers to the bits less significant than $M^H$. For example, $M = 1.\{x\ldots\}$, $M^H = 1.\{x\}_k$, $M^L = 0.\{0\}_k\{x\ldots\}$.

We derive an upper bound on $M_C^H - M'_C$ in three steps:

**Step 1**: Analyze the relationship between $M_{\hat{C}}^H$ and $(M'_{\hat{C}})^H$.

**Step 2**: Analyze (a) how rounding in the FPU affects $M_C^H$ and (b) how rounding in the checker affects $M'_C$.

**Step 3**: Analyze the relationship between $M_C^H$ and $M'_C$.

**Step 1**

Given $E_A + E_B = E_C$ in Eqn. (2), the unrounded product of the mantissas is:

$$M_A M_B = M_{\hat{C}} = M_{\hat{C}}^H + M_{\hat{C}}^L$$

In Eqn. (4), the checker computes the unrounded product:

$$M_A^H M_B^H = M'_{\hat{C}} = (M'_{\hat{C}})^H + (M'_{\hat{C}})^L$$

From Eqn. (2), we can write:

$$M_{\hat{C}} = M_A M_B = M'_{\hat{C}} + M^*$$

where $M^*$ was defined to be $M_A^H M_B^L + M_A^L M_B^H + M_A^L M_B^L$ and $M'_{\hat{C}} = M_A^H M_B^H$. Again, note that $M^*$ can have more than *k* bits after the binary point. Splitting the mantissas into High and Low on both sides, we can write the previous expression as

$$M_{\hat{C}}^H + M_{\hat{C}}^L = (M'_{\hat{C}})^H + (M'_{\hat{C}})^L + (M^*)^H + (M^*)^L \qquad \text{Eqn. (33)}$$

By Axiom 3, a Low term is always smaller than $2^{-k}$, so $(M'_{\hat{C}})^L + (M^*)^L$ is in the range $[0, 2 \times 2^{-k})$. We sub-divide this range into two sub-ranges for purposes of our analysis:

Range 1: $0 \le (M'_{\hat{C}})^L + (M^*)^L < 2^{-k}$

Range 2: $2^{-k} \le (M'_{\hat{C}})^L + (M^*)^L < 2 \times 2^{-k}$

In Range 1, the sum of $(M'_{\hat{C}})^L + (M^*)^L$ does not cause a carry into the sum of the High terms $(M'_{\hat{C}})^H + (M^*)^H$ in Eqn. (33). So:

$$\begin{cases} M_{\hat{C}}^H = \left(M'_{\hat{C}}\right)^H + (M^*)^H \\ M_{\hat{C}}^L = \left(M'_{\hat{C}}\right)^L + (M^*)^L \end{cases} \qquad \text{Case I}$$

In Range 2, the sum of $\left(M'_{\hat{C}}\right)^L + (M^*)^L$ causes a carry into the sum of the High terms. So $M_{\hat{C}}^H = \left(M'_{\hat{C}}\right)^H + (M^*)^H + 2^{-k}$. And

$$M_{\hat{C}}^L = M_{\hat{C}} - M_{\hat{C}}^H = \left(M'_{\hat{C}} + M^*\right) - \left(\left(M'_{\hat{C}}\right)^H + (M^*)^H + 2^{-k}\right)$$
$$= \left(M'_{\hat{C}}\right)^L + (M^*)^L - 2^{-k}$$

Thus,

$$\begin{cases} M_{\hat{C}}^H = \left(M'_{\hat{C}}\right)^H + (M^*)^H + 2^{-k} \\ M_{\hat{C}}^L = \left(M'_{\hat{C}}\right)^L + (M^*)^L - 2^{-k} \end{cases} \qquad \text{Case II}$$

**Step 2.a**

For the FPU, $M_{\hat{C}} = M_C + \delta_C$. So

$$M_{\hat{C}} = M_{\hat{C}}^H + M_{\hat{C}}^L = M_C^H + M_C^L + \delta_C$$

The only way in which $M_{\hat{C}}^H \neq M_C^H$ is if, after rounding, the carried one at the 23rd bit propagates to the most significant *k* bits. For this to happen, the (*k*+1)th to the 23rd bits after the binary point in $M_{\hat{C}}^L$ must all equal 1. To produce the carried one, the bits after the 23rd bit of $M_{\hat{C}}$ must be larger than $2^{-24}$, so $M_{\hat{C}}$ rounds up to get $M_C$. Thus, $M_{\hat{C}}^H \neq M_C^H$ when $M_{\hat{C}}^L > 0.\{0\}_k\{1\}_{23-k}\{1\} = 2^{-k} - 2^{-24}$. After rounding, all bits in $M_C^L$ must be zero because they all sum with the carried one. For purposes of analysis, we consider two ranges of values for $M_{\hat{C}}^L$:

Range 1: $\quad 0 < M_{\hat{C}}^L < 2^{-k} - 2^{-24}$

Range 2: $\quad 2^{-k} - 2^{-24} < M_{\hat{C}}^L < 2^{-k} \text{ and } M_C^L = 0$

In Range 1, $M_C^H$ is equivalent to $M_{\hat{C}}^H$ after rounding:

$$\begin{cases} M_{\hat{C}}^H = M_C^H \\ M_{\hat{C}}^L = M_C^L + \delta_C \end{cases} \qquad \text{Case III}$$

In Range 2, $M_{\hat{C}}$ rounds up to $M_C^H$ and is $2^{-k}$ larger than $M_{\hat{C}}^H$ after rounding:

$$\begin{cases} M_{\hat{C}}^H = M_C^H - 2^{-k} \\ M_{\hat{C}}^L = M_C^L + \delta_C + 2^{-k} \\ \quad = \delta_C + 2^{-k} \end{cases} \qquad \text{Case IV}$$

**Step 2.b**

For the checker,

$$M'_{\hat{C}} = \left(M'_{\hat{C}}\right)^H + \left(M'_{\hat{C}}\right)^L = M'_C + \delta_C$$

Notice there is no $(M'_C)^L$ term, because $M'_C$ is *produced* by the checker and has only *k* bits after the binary point. For analysis, we consider two ranges of values for $\left(M'_{\hat{C}}\right)^L$:

Range 1: $\quad 0 \leq \left(M'_{\hat{C}}\right)^L < 2^{-k-1}$

Range 2: $\quad 2^{-k-1} \leq \left(M'_{\hat{C}}\right)^L < 2^{-k}$

In Range 1, $M'_{\hat{C}}$ rounds down to get $M'_C$. $M'_C$ equals $\left(M'_{\hat{C}}\right)^H$ after rounding.

$$\begin{cases} \left(M'_{\hat{C}}\right)^H = M'_C \\ \left(M'_{\hat{C}}\right)^L = \delta'_C \end{cases} \qquad \text{Case V}$$

In Range 2, $M'_{\hat{C}}$ rounds up to get $M'_C$. $M'_C$ is $2^{-k}$ greater than $\left(M'_{\hat{C}}\right)^H$ after rounding.

$$\begin{cases} \left(M'_{\hat{C}}\right)^H = M'_C - 2^{-k} \\ \left(M'_{\hat{C}}\right)^L = \delta'_C + 2^{-k} \end{cases} \qquad \text{Case VI}$$

**Step 3**

We now have to consider all possible combinations of cases from Steps 1, 2.a, and 2.b. Because there are two cases in each step, there are $2^3$=8 theoretical combinations, but some combinations are impossible.

- Case I-Case III-Case V: Plug in the equations for the high bits in Case III and Case V into Case I.
$$M_C^H = M'_C + (M^*)^H$$
Rearrange the equation to get:
$$M_C^H - M'_C = (M^*)^H < 4 \times 2^{-k}$$

- Case I-Case IV-Case V (impossible): $M_{\hat{C}}^L = \left(M'_{\hat{C}}\right)^L + (M^*)^L$ from Case I. $M_{\hat{C}}^L$ is in the range $(2^{-k} - 2^{-24}, 2^{-k})$, and $\left(M'_{\hat{C}}\right)^L$ is in the range $[0, 2^{-k-1})$. Then the boundary condition for $(M^*)^L$ is
$$2^{-k} - 2^{-24} < (M^*)^L < 2^{-k-1}$$
Since the upper bound is smaller than the lower bound, because *k*<24, there is no $(M^*)^L$ that satisfies the requirements of Case I-Case IV-Case V, i.e. the case does not exist.

- Case I-Case III-Case VI: Plug Case III and Case VI into Case I:
$$M_C^H = M'_C + (M^*)^H - 2^{-k}$$
$$M_C^H - M'_C = (M^*)^H - 2^{-k} < 3 \times 2^{-k}$$

- Case I-Case IV-Case VI: Plug in Case Case IV and Case VI into Case I:
$$M_C^H - 2^{-k} = M'_C + (M^*)^H - 2^{-k}$$
$$M_C^H - M'_C = (M^*)^H < 4 \times 2^{-k}$$

- Case II-Case III-Case V (impossible): From Case II, $\left(M'_{\hat{C}}\right)^L + (M^*)^L = M_{\hat{C}}^L + 2^{-k}$. $\left(M'_{\hat{C}}\right)^L + (M^*)^L$ is in the range $(2^{-k}, 2 \times 2^{-k} - 2^{-24})$, and $\left(M'_{\hat{C}}\right)^L$ is in the range $[0, 2^{-k-1})$. Then the boundary condition of $(M^*)^L$ is
$$2^{-k} < (M^*)^L < 1.5 \times 2^{-k} - 2^{-24}$$
which is invalid because $(M^*)^L < 2^{-k}$. Therefore, Case II-Case III-Case V is impossible.

- Case II-Case III-Case VI: Plug Case III and Case VI into Case II:
$$M_C^H - 2^{-k} = M'_C + (M^*)^H - 2^{-k}$$
$$M_C^H - M'_C = (M^*)^H < 4 \times 2^{-k}$$

- Case II-Case III-Case V (impossible): $\left(M'_{\hat{C}}\right)^L + (M^*)^L = M_{\hat{C}}^L + 2^{-k}$ from Case II. $\left(M'_{\hat{C}}\right)^L + (M^*)^L$ is in the range $(2 \times 2^{-k} - 2^{-24}, 2 \times 2^{-k})$. $\left(M'_{\hat{C}}\right)^L$ is in the range $[0, 2^{-k-1})$. Then
$$2 \times 2^{-k} - 2^{-24} < (M^*)^L < 1.5 \times 2^{-k}$$
Because the upper bound is smaller than the lower bound, this case is impossible.

- Case II-Case IV-Case VI (impossible): As in Case II-Case III-Case V, $+(M^*)^L$ is in the range $(2 \times 2^{-k} - 2^{-24}, 2 \times 2^{-k})$. $\left(M'_{\hat{C}}\right)^L$ is in the range $[2^{-k-1}, 2^{-k})$. Thus
$$1.5 \times 2^{-k} - 2^{-24} < (M^*)^L < 1.5 \times 2^{-k}$$
Because $(M^*)^L < 2^{-k}$, this case is impossible.

**Conclusion** Overall, $M_C^H - M_C' < 4 \times 2^{-k}$. Therefore $Diff = (M_C^H - M_C') \times 2^{-k}$ in Section 5.1 is in [-1,3].